\documentclass[twoside,draft,reqno]{birkart}

\usepackage{amssymb}

\let\cal\mathcal
\let\frak\mathfrak
\let\Bbb\mathbb
\def\HS{{\rm HS}}
\def\hsp{{\rm hsp}}
\def\top{{\rm top}}
\def\Var{{\rm Var}}
\def\top{{\rm top}}
\def\Gal{{\rm Gal}}

\newtheorem{theorem}{Theorem}[subsection]

\newtheorem{corollary}[theorem]{Corollary}

\newtheorem{conjecture}[theorem]{Conjecture}

\theoremstyle{definition}
\newtheorem{definition}[theorem]{Definition}

\newtheorem{notation}[theorem]{Notation}

\theoremstyle{plain}

\begin{document}
\setcounter{page}{1}

\title[Geometry on arc spaces of algebraic varieties]
{Geometry on arc spaces of algebraic varieties}

\author{Jan Denef}
\address{University of Leuven, Department of Mathematics,
Celestijnenlaan 200B, 3001 Leu\-ven, Bel\-gium }
\email{ Jan.Denef@wis.kuleuven.ac.be}
\urladdr{http://www.wis.kuleuven.ac.be/wis/algebra/denef.html}

\author{Fran\c cois Loeser}

\address{Centre de Math{\'e}matiques,
Ecole Polytechnique,
F-91128 Palaiseau
(UMR 7640 du CNRS), {\rm and}
Institut de Math{\'e}matiques,
Universit{\'e} P. et M. Curie, Case 82,
4 place Jussieu,
F-75252 Paris Cedex 05
(UMR 7596 du CNRS)}
\email{loeser@math.polytechnique.fr}
\urladdr{http://math.polytechnique.fr/cmat/loeser/loeser.html}

%\author[J.~Denef and F. Loeser]{J. Denef and F. Loeser}

%\address{%%
%Jan Denef,\br
%Department of Mathematics,\br
%Katholieke Universiteit Leuven,\br
%Celestijnenlaan 200B \br
%3001 Leuven (Heverlee) BELGIUM \br
%{\it E-mail address :} {\bf jan.denef@wis.kuleuven.ac.be} \br
%http : }

% \email{jan.denef@wis.kuleuven.ac.be}

%\address{%%
%Fran\c cois Loeser,\br
%Centre de Math{\'e}matiques,\br
%Ecole Polytechnique,\br
%F-91128 Palaiseau (UMR 7640 du CNRS),\br
%Institut de Math{\'e}matiques,\br
%Universit{\'e} P. et M. Curie, Case 82, 4 place Jussieu,\br
%F-75252 Paris Cedex 05 (UMR 7596 du CNRS),\br
%{\it E-mail address :} {\bf loeser@math.polytechnique.fr} \br
%http : }

\begin{abstract}
This paper is a survey on arc spaces, a recent topic in algebraic
geometry and singularity theory. The geometry of
the arc space of an algebraic variety yields several new geometric
invariants and brings new light to some classical invariants.

\end{abstract}

\maketitle

\section{Introduction}\label{sec1}For an algebraic variety $X$ over the field
$\Bbb C$ of complex numbers, one considers the arc space ${\cal L}(X)$, whose points are the ${\Bbb C}[[t]]$-rational points on $X$, and the
truncated arc spaces ${\cal L}_n (X)$, whose points are the ${\Bbb
C}[[t]]/t^{n +1}$-rational points on $X$. The geometry of these spaces
yields several new geometric invariants of $X$ and brings new
light to some classical invariants. For example, Denef and
Loeser \cite{DeLo2} showed that the Hodge spectrum of a critical point
of a polynomial can be expressed in terms of geometry on arc
spaces, yielding a new proof and a generalization \cite{DeLo4}
of the Thom-Sebastiani Theorem for the Hodge spectrum due
to Varchenko \cite{Va}
and Saito \cite{Sa3}, 
\cite{Sa4}.
In a
different direction, Batyrev \cite{Ba3}
used arc spaces to prove a
conjecture of Reid \cite{Re} on quotient singularities (the McKay
correspondence), and to construct his stringy Hodge numbers \cite{Ba2}
appearing in mirror symmetry. All these developments are based on
Kontsevich's construction \cite{Ko} of a measure on the arc space
${\cal L}(X)$, the motivic measure, which is an analogue of
the $p$-adic measure on a $p$-adic variety.

In section \ref{sec2} we define the arc spaces ${\cal L}(X)$ and
${\cal L}_n(X)$ of an algebraic variety over any field $k$ of characteristic
zero. The first question that appears is how ${\cal L}_n(X)$ and
$\pi_n({\cal L}_n(X))$ change with $n$, where $\pi_n$ denotes the
truncation map from ${\cal L}(X)$ to ${\cal L}_n(X)$. As a partial
answer to this question we will see in 2.2.1 that the power
series
 $$J(T,\chi) := \sum_{n \geq 0} \, \chi
({\cal L}_n(X)) \,T^n, \quad P(T,\chi) = \sum_{n \geq 0} \,
\chi(\pi_n({\cal L}(X)))\,T^n$$ are rational (i.e. a quotient of two polynomials), for
any reasonable generalized Euler characteristic $\chi$. This is a
direct consequence of results of Denef and Loeser \cite{DeLo3}. Instead
of working with particular generalized Euler characteristics, such
as the topological Euler characteristic, the Hodge polynomial or
the Hodge characteristic, it is more general to work with the
universal Euler characteristic which associates to any algebraic
variety $X$ over $k$ its class $[X]$ in the Grothendieck group
$K_0(\text{Var}_k)$ of algebraic varieties over $k$. This is the
abelian group generated by symbols $[X]$, for $X$ a variety over
$k$, with the relations $[X] = [Y]$ if $X$ and $Y$ are isomorphic,
and $[X] = [Y] + [X \setminus Y]$ if $Y$ is Zariski closed in $X$.
There is a natural ring structure on $K_0(\text{Var}_k)$, the
product of $[X]$ and $[Y]$ being equal to $[X \times Y]$. We denote
by $\mathcal{M}_k$ the ring obtained from $K_0(\text{Var}_k)$ by
inverting the class of $\Bbb A_k^1$. The above rationality result
applied to the universal Euler characteristic says that the power
series $$J(T) := \sum_{n \geq 0} \, [{\cal L}_n(X)] \,T^n, \quad P(T) :=
\sum_{n \geq 0}
\,
[\pi_n({\cal L}(X))] \,T^n$$ in $\mathcal{M}_k [[T]]$ are rational.

Power series like $J(T)$ and $P(T$), with coefficients in
$\mathcal{M}_k[[T]]$, are called ``motivic'', because they
specialize to power series over the Grothendieck group $K_0
({\rm Mot}_k)$ of the category of Chow motives over $k$. Actually in
several of our papers on arcs we work over $K_0 ({\rm Mot}_k)$ instead
of over $\mathcal{M}_k[[T]]$.

In section \ref{sec3} we introduce the motivic zeta function $Z(T)$
associated to a morphism $f$ from an nonsingular algebraic variety
$X$ to the affine line, cf. \cite{DeLo7}. A naive version of it is the
power series over $\mathcal{M}_k$ defined by $$Z^{\text{naive}}(T)
:= \sum_{n \geq 1} \, [\frak X_n ]\, [\Bbb A_k^1 ]^{-nd} \, T^n.$$
Here
$\frak X_n$ denotes the set of arcs $\varphi$ in ${\cal L}(X)$ with
$f(\varphi)$ a power series of order $n$. The motivic zeta function
of $f$ contains a wealth of geometric information about $f$. For
example the Hodge spectrum of any critical point of $f$ can be
expressed in terms of $\lim_{T \rightarrow \infty} Z(T)$. This
limit is a well defined element of $\mathcal{M}_k$, and can be
considered as the ``virtual motivic incarnation'' of the Milnor
fibers of $f$. All this is explained in section 3.5. In section
3.4 we also show that the topological zeta functions of Denef and
Loeser \cite{DeLo1} can be expressed in terms of the motivic zeta
function.

We explain in section \ref{sec4} the notion of motivic integration on
${\cal L}(X)$, due to Kontsevich \cite{Ko}, and further developed
by Batyrev \cite{Ba2}, \cite{Ba3}, Denef and
Loeser [DeLo2-7], and Looijenga \cite{Loo}.
This notion plays a key role in the present paper.
Kontsevich used it to prove that two birationally
equivalent Calabi-Yau manifolds have the same Hodge numbers. This
result, together with
some other direct applications of motivic integration, is
discussed in section 4.4.

One of the most striking applications of arc spaces and motivic
integration is Batyrev's proof \cite{Ba3} of the conjecture of Reid on
the generalized McKay correspondence. We will not treat this
material in the present paper, but refer to the Bourbaki report of
Reid \cite{Re}, see also \cite{DeLo5} and \cite{Loo}.

In section \ref{sec5} we explain how the relation between the Hodge
spectrum and the motivic zeta function yields a new proof of
Varchenko's and
Saito's Thom-Sebastiani Theorem which expresses the Hodge spectrum
of $ f(x)+g(y)$ in terms of the Hodge spectra of $f(x)$ and
$g(y)$. Our method \cite{DeLo4} actually yields a much stronger result
which relates the ``virtual motivic incarnations'' of the Milnor
fibers of these three functions. In our paper \cite{DeLo4} we obtained
this result only at the level of Chow motives, and is was
Looijenga \cite{Loo} who showed how to work at the level  of the
Grothendieck ring of algebraic varieties.

Finally in section \ref{sec6} we briefly discuss the connections with the
$p$-adic case. Considering  motivic integration as an analogue of
$p$-adic integration, several arithmetical results in the $p$-adic
case find their natural counterpart in complex geometry and in
the theory of motives.

In the present paper we have avoided to work with Chow motives (with
one important exception in section \ref{sec6}). Indeed, the more
recent material in \cite{DeLo7} and \cite{Loo} shows that this is
possible except for the functional equation in section 3 of 
\cite{DeLo2} and the results in \cite{DeLo6}.

\section{The arc space of a variety}\label{sec2}
We fix a base field $k$ of characteristic zero.  The reader may choose to only
consider the case where $k$ is the field $\Bbb C$ of complex numbers.  Let $X$
be an algebraic variety over $k$, not necessarily irreducible, i.e. $X$ is a
reduced separated scheme of finite type over $k$.

\subsection{The arc space of $X$}
For each natural number $n$ we consider the {\it space ${\cal L}_n(X)$ of arcs
modulo
$t^{n+1}$ on $X$}. This is an algebraic variety over $k$, whose $K$-rational
points, for any field $K$ containing $k$, are the $K[t]/t^{n+1}K[t]$-rational
points of $X$.  For example when $X$ is an affine variety with equations $f_i
(\vec x) = 0$, $i= 1, \cdots, m$, $\vec x = (x_1,\cdots, x_n)$, then ${\cal L}_n(X)$
is given by the equations, in the variables $\vec a_0, \cdots, \vec a_n$,
expressing that $f_i (\vec a_0 + \vec a_1t + \cdots + \vec a_n t^n) \equiv 0
\mod t^{n+1}, i = 1,\cdots, m$.  

Taking the projective limit of these algebraic varieties ${\cal L}_n(X)$ we
obtain
the {\it arc space ${\cal L}(X)$ of} $X$, which is a reduced separated scheme
over
$k$.  In general, ${\cal L}(X)$ is not of finite type over $k$ 
(i.e. ${\cal L}(X)$ is an ``algebraic
variety of infinite dimension'').  The $K$-rational points of ${\cal L} (X)$ are
the
$K[[t]]$-rational points of $X$.  These are called {\it $K$-arcs on $X$}.  For
example when $X$ is an affine complex variety with equations $f_i(\vec x) = 0,
i = 1, \cdots m, \vec x = (x_1,\cdots,x_n)$, then the $\Bbb C$-rational points
of
${\cal L}(X)$ are the sequences $(\vec a_0, \vec a_1, \vec a_2, \cdots) \in
(\Bbb C^n)^{\Bbb N}$ satisfying $f_i (\vec a_0 + \vec a_1 t + \vec a_2t^2 +
\cdots )  = 0$,  for $i = 1,\cdots,m$. 
For any $n$, and for $m > n$, we have natural morphisms
$$\pi_n : {\cal L}(X) \rightarrow {\cal L}_n(X) \ \ {\rm and } \ \ \pi^m_n :
{\cal L}_m(X) \rightarrow {\cal L}_n(X),$$
obtained by truncation.  Note that ${\cal L}_0(X) = X$ and that ${\cal L}_1(X)$
is the tangent bundle of $X$.  For any arc $\gamma$ on $X$ (i.e. a $K$-arc for
some field $K$ containing $k$), we call $\pi_0(\gamma)$ {\it the origin of the
arc
$\gamma$}.

By a theorem of Greenberg \cite{Gr}, given an algebraic variety $X$ over $k$, there
exists a number $c > 0$, such that for any $n$ and for any field $K$ containing
$k$ we have
$$\pi_n ({\cal L}(X)(K))  = \pi^{cn}_n({\cal L}_{cn} (X)(K)),$$
writing $Y(K)$ to denote the set of $K$-rational points on any variety
$Y$ over $k$.  This implies that $\pi_n({\cal L}(X))$ is a constructible subset
of the algebraic variety ${\cal L}_n(X)$.  If $X$ is smooth, then we can take
$c=1$ and $\pi_n$ is surjective.  Moreover, in that case, $\pi^m_n$ is a
locally trivial fibration with fiber $\Bbb A^{(m-n)\dim X}_{k}$. Here $\Bbb A^d_k$
denotes the affine space of dimension $d$ over $k$.

Probably Nash \cite{Na} was the first to study arc spaces in a systematic way (his
paper was written in 1968, but published only recently). 
For a singular point
$P$ on $X$, he considered the space ${\cal L}_{\{P\}}(X) := \pi^{-1}_0(P)$ of
arcs on $X$ with origin $P$,
and its subspace
${\cal N}_{\{P\}}(X)$ of arcs with origin $P$ which are
not contained in the
singular locus of $X$.
He proved that the number of irreducible
components of the Zariski closure of $\pi_n({\cal N}_{\{P\}}(X))$ stabilizes
for $n$ big enough, and associated to each irreducible component of this
closure (when $n \gg 0$), in a canonical and injective way, an irreducible
component of the preimage of $P$ in any resolution of singularities of $X$.
Moreover he conjectured that any irreducible component of the preimage of $P$
in a given resolution of $X$, which ``appears'' in all resolutions of $X$, is
obtained in that way. For other results concerning arc spaces we 
refer to \cite{Le},
\cite{Hi} and \cite{LeRe}. 

\subsection{How do ${\cal L}_n(X)$ and $\pi_n({\cal L}(X))$ change with $n$ ?}
The work of  Nash \cite{Na} is the first result towards the question of how the
geometry of $\pi_n({\cal L}_{\{ P \}}(X))$ changes with $n$.  Recently we
investigated how the topological Euler characteristic $\chi_{\top}$ (case
$k = \Bbb C$), and
generalized Euler characteristics of the spaces ${\cal L}_n(X),
\pi_n({\cal L}(X)), \pi_n ({\cal L}_{\{P\}}(X))$ change with $n$.

With a {\it generalized Euler characteristic} on the category $\Var_k$ of
algebraic
varieties over $k$, we mean a map $\chi$ from $\Var_k$ to some commutative ring
$R$ such that $\chi(X) = \chi (Y)$ when $X \cong Y, \ \chi(X) = \chi(Y) +
\chi(X
\setminus Y)$ when $Y$ is a Zariski closed subvariety of $X$, and $\chi(X
\times Y) = \chi(X) \cdot \chi(Y)$.
\smallskip
Clearly $\chi = \chi_{\top}$ satisfies the above requirements with $R = \Bbb
Z$, when $k = \Bbb C$.  Another example of a generalized Euler characteristic,
when $k = \Bbb C$, is given by $\chi_{hp} :
\Var_{\Bbb C} \rightarrow \Bbb
Z[u,v] : X \mapsto \sum_{i,p,q} (-1)^i h^i_{p,q}u^pv^q$, where $h^i_{p,q}$ is
the
dimension of the $(p,q)$-component of the mixed Hodge structure on $H^i_c(X,
\Bbb C)$.  One calls $\chi_{hp}(X)$ the {\it Hodge polynomial of} $X$.  For
example, $\chi_{hp}(\Bbb A^1_{\Bbb C}) = uv$.  We refer to 3.1.2 for the {\it
Hodge characteristic} $\chi_h$ which takes values in the Grothendieck group
$K_0(\HS)$ of the abelian category of Hodge structures.  There are many other
examples of generalized Euler characteristics.  For example, to mention an
exotic one, when $k = \Bbb Q$, there is {\it the conductor} $c(X)$ of $X$,
which yields a generalized Euler characteristic $c :
{Var}_{\Bbb Q} \rightarrow
\Bbb Q_{> 0} : X \mapsto c(X) := \prod_i (c_i)^{{(-1)}^{i+1}}$, where $c_i$
denotes the conductor of the $\ell$-adic representation of
$\Gal(\bar {\Bbb
Q},
\Bbb Q)$ on
the {\'e}tale cohomology $H^i_c(X_{\bar{\Bbb Q}}, \Bbb Q_\ell)$, where $\ell$
is a
fixed prime
and $\bar{\Bbb Q}$ an algebraic closure of $\Bbb Q$; see e.g. \cite{Se}.
Here $\Bbb
Q_{>
0}$ is
the multiplicative group of positive rational numbers with ring structure
inherited by $Q_{> 0} \cong \oplus_{p \, {\text{prime}}} \Bbb Z : x \mapsto
({\rm ord}_p x)_p$.  For example, when $X$ is an elliptic curve, $c(X)$ is the
usual conductor of $X$, related to the primes at which $X$ has had reduction.

\begin{theorem}[\cite{DeLo3}]\label{2.2.1}
Let $\chi : \Var_k \rightarrow R$ be a generalized Euler
characteristic,
and suppose that $\chi(\Bbb A^1_k)$ is not a zero divisor in $R$.  Then the
power series
$$J(T,\chi) := \sum_{n \geq 0} \, \chi ({\cal L}_n(X)) \,
T^n, \quad P(T,\chi)
= \sum_{n \geq 0} \,
\chi(\pi_n({\cal L}(X))) \,T^n$$
are rational (i.e. a quotient of two polynomials).  Actually the denominators
are products of polynomials of the form $1 - (\chi(\Bbb A^1_k))^a T^b$, with $b
\in \Bbb N \setminus \{ 0 \}, a \in \Bbb Z$.
\end{theorem}

In 2.3 below, we will construct the universal Euler characteristic
${\Var}_k \rightarrow K_0(\Var_k) : X \mapsto [X]$, where
$K_0(\Var_k)$ denotes the Grothendieck group of varieties over $k$ (see
2.3).  Any generalized Euler characteristic on $\Var_k$ factorizes over this
universal one.  So it suffices to prove Theorem \ref{2.2.1} with $R$ the ring
$\mathcal{M}_k$
obtained from $K_0 (\Var_k)$ by inverting $[\Bbb A^1_k]$.  The rationality of
$J(T,\chi)$ follows from the material we will discuss in section 3, when $X$ is
an
affine hypersurface, see 3.3.3.  The proof of the rationality of $P(T,\chi)$ is
more complicated and uses a theorem of Pas \cite{Pa} on quantifier elimination for
power series rings.  If $f : X \rightarrow Y$ is a morphism of algebraic
varieties, then the image $f(A)$ in $\mathcal{L}(Y)$ of a constructible subset
$A$ (cf. 4.1 below) of $\mathcal{L}(X)$ is generally not constructible.  The
theorem of Pas implies that $f(A)$ still has a ``simple description'' and is in
fact what is called a semi-algebraic subset of $\mathcal{L}(Y)$, cf. 
\cite{DeLo3}.
This
plays a key role in the proof of the rationality of $P(T,\chi)$.  For 
a survey on applications of quantifier elimination results for valued 
fields, see \cite{De4}.

\subsection{Grothendieck groups of varieties}

Let $S$ be an algebraic variety over $k$.  By an $S$-variety we mean a
variety $X$ together with a morphism $X \rightarrow S$.  The $S$-varieties from
a category denoted by $\Var_S$, the arrows are the morphisms that commute with
the morphisms to $S$.

We denote by $K_0 (\Var_S)$ {\it the Grothendieck group of
$S$-varieties}. It is
an abelian group generated by symbols $[X]$, for $X$ an $S$-variety, with the
relations $[X] = [Y]$ if $X$ and $Y$ are isomorphic in
$\Var_S$, and $[X] = [Y]
+ [X \setminus Y]$ if $Y$ is Zariski closed in $X$.  There is a natural ring
structure on $K_0 (\Var_S)$, the product of $[X]$ and $[Y]$ being equal to $[X
\times_S Y]$.  Sometimes we will also write $[X/S]$ instead of $[X]$, to
emphasize the role of $S$.  We write $\Bbb L$ to denote the class of $\Bbb
A^1_k \times S$ in $K_0(\Var_S)$, where the morphism from $\Bbb A^1_k
\times S$ to $S$ is the natural projection.  We denote by $\mathcal{M}_S$ the
ring
obtained from $K_0(\Var_S)$ by inverting $\Bbb L$.  When $A$ is a
constructible subset of some $S$-variety, we define $[A/S]$ in the obvious way,
writing $A$ as a disjoint union of a finite number of
locally closed subvarieties $A_i$.  Indeed
$[A/S] := \sum_i [A_i/S]$ does not depend on the choice of the
subvarieties
$A_i$. 

When $S$ consists of only one geometric point, i.e. $S = {\rm Spec} (k)$, then
we will write $K_0(\Var_k)$ instead of $K_0(\Var_S)$ (to denote the
Grothendieck group of algebraic varieties over $k$), and $\mathcal{M}_k$
instead of
$\mathcal{M}_S$. 
Clearly the map $\Var_k \rightarrow K_0(\Var_k)$ is the universal
generalized Euler characteristic, in the sense that any generalized Euler
characteristic on $\Var_k$ factors through it.  

In our papers [DeLo2-7], we always work with $K_0(\Var_k)$, but
recently E. Looijenga, in his Bourbaki talk \cite{Loo}, introduced the relative
Grothendieck ring $K_0(\Var_S)$, stating some of our results in a
stronger form.  For example, considering ${\cal L}_n(X)$ as an $X$-variety
through the morphism $\pi_0^n$, our proof of Theorem \ref{2.2.1} actually yields the
slightly stronger

\begin{theorem}\label{2.3.1}
Let $X$ be a variety over $k$, then the power series 
$$J(T) := \sum_{n \geq 0} \, [{\cal L}_n(X)/X] \, T^n, \quad P(T) := \sum_{n
  \geq 0} \, [\pi_n({\cal L}(X))/X] \, T^n$$
in $\mathcal{M}_X[[T]]$ are rational, with denominator a product of polynomials
of
the form $1 - \Bbb L^a T^b$, with $b \in \Bbb N \setminus \{ 0 \}, a \in \Bbb
Z$.
\end{theorem}

\subsection{Equivariant Grothendieck groups}
We need some technical preparation in order to take care of the monodromy
actions in the next section.

For any positive integer $n$, let $\mu_n$ be the group of all $n$-th roots of
unity (in some fixed algebraic closure of $k$).  Note that $\mu_n$ is actually
an algebraic variety over $k$, namely ${\rm Spec} (k[x]/(x^n-1))$.
The $\mu_n$ form a projective system, with respect to the maps $\mu_{nd}
\rightarrow \mu_n : x \mapsto x^d$.  We denote by $\hat \mu$ the projective
limit of the $\mu_n$.  Note that the group $\hat \mu$ is not an algebraic
variety. It is called a pro-variety.

Let $X$ be an $S$-variety.  A {\it good $\mu_n$-action on $X$} is a group
action
$\mu_n \times X \rightarrow X$ which is a morphism of $S$-varieties, such that
each orbit is contained in an affine subvariety of $X$.  This last condition is
automatically satisfied when $X$ is a quasi projective variety.  A {\it good
$\hat
\mu$-action on $X$} is an action of $\hat \mu$ on $X$ which factors through a
good $\mu_n$-action, for some $n$.

The {\it monodromic Grothendieck group} $K^{\hat \mu}_0
(\Var_S)$ is defined as
the
abelian group generated by symbols $[X, \hat \mu]$ (also denoted by $[X/S,\hat
\mu]$, or simply $[X]$), for $X$ an $S$-variety with good $\hat \mu$-action,
with the relations $[X, \hat \mu]  = [Y,\hat \mu]$ if $X$ and $Y$ are
isomorphic as $S$-varieties with $\hat \mu$ -ction, and $[X,\hat \mu] = [Y,\hat
\mu] + [X \setminus Y, \hat \mu]$ if $Y$ is Zariski closed in $X$ with the
$\hat \mu$-action on $Y$ induced by the one on $X$, and moreover $[X \times
V,\hat \mu] = [X \times \Bbb A^n_k, \hat \mu]$ where $V$ is the $n$-dimensional
affine space over $k$ with any good $\hat \mu$-action, and $\Bbb A^n_k$ is
taken with the trivial $\hat \mu$-action.  There is a natural ring structure on
$K^{\hat \mu}_0 (\Var_S)$, the product being induced by the fiber product over
$S$.  We write $\Bbb L$ to denote the class in $K_0^{\hat \mu}
(\Var_S)$ of
$\Bbb A^1_k \times S$ with the trivial $\hat \mu$-action.

We denote by $\mathcal{M}^{\hat \mu}_S$ the ring obtained from $K_0^{\hat
\mu}(\Var_S)$ by inverting $\Bbb L$.  When $A$ is a constructible subset
of $X$ which is stable under the $\hat \mu$-action, then we define $[A,\hat
\mu]$ in the obvious way.  When $S$ consists of only one geometric point, i.e.
$S = {\rm Spec} (k)$, then we will write $K^{\hat \mu}_0(\Var_k)$
instead of $K^{\hat \mu}_0(\Var_S)$. The group $K^{\hat \mu}_0
(\Var_k)$
was first introduced in \cite{DeLo7}. 

Note that for any $s \in S(k)$ we have natural maps
$K^{\hat \mu}_0 (\Var_S) \rightarrow K^{\hat \mu}_0 (\Var_k)$
and $\mathcal{M}^{\hat \mu}_S \rightarrow \mathcal{M}^{\hat \mu}_k$
given by
$[X,\hat
\mu] \rightarrow [X_s,\hat \mu]$, where $X_s$ denotes the fiber at $s$ of $X
\rightarrow S$.

Although $\mathcal{M}^{\hat \mu}_k$ is a very complicated ring, there are many
interesting morphisms from it to simpler rings.  For example when $k = \Bbb C$,
for any character $\alpha$ of $\hat \mu$ (i.e. a group homomorphism $\alpha :
\hat \mu \rightarrow \Bbb C^\times$), there is a natural ring homomorphism
$$\chi_{\top}(-,\alpha) : \mathcal{M}^{\hat \mu}_k \longrightarrow
\Bbb Z :
\quad
X
\longmapsto \sum_{q \geq 0} (-1)^q {\rm dim } H^q(X,\Bbb C)_\alpha,$$
where $H^\ast(X, \Bbb C)_\alpha$ is the part of $H^\ast(X,\Bbb C)$ on which
$\hat \mu$ acts by multiplication by $\alpha$.

\section{The motivic zeta function of a regular function}\label{sec3}
Let $X$ be a nonsingular irreducible algebraic variety over $k$ of dimension
$d$ and $f : X \rightarrow \Bbb A^1_k$ a non constant
morphism.  In this section we introduce several new invariants of $f$.  These
are constructed using arc spaces.  We first recall in 3.1 some classical
invariants associated to $f$.  In what follows we denote by $X_0$ the {\it
locus of
$f=0$} in $X$.

\subsection{The monodromy zeta function and the Hodge spectrum}
In this subsection 3.1 we suppose that $k = \Bbb C$.  Let $x$ be a point of
$X_0 =
f^{-1}(0)$.  We fix a smooth metric on $X$.

\subsubsection{Monodromy.}\label{3.1.1}
We set $X^\times_{\epsilon,\eta} := B(x,\epsilon) \cap
f^{-1}(D^\times_\eta)$, with $B(x,\epsilon)$ the open ball of radius $\epsilon$
centered at $x$ and $D^\times_\eta := D_\eta \setminus \{ 0 \}$, with $D_\eta$
the open disk of radius $\eta$ centered at $0$.  For $0 < \eta \ll \epsilon \ll
1$,
the restriction of $f$ to $X^\times_{\epsilon,\eta}$ is a locally trivial
fibration, called the {\it Milnor fibration,} onto $D^\times_\eta$ with fiber
$F_x$, {\it the Milnor fiber at} $x$.  The action of a characteristic
homeomorphism of this fibration on cohomology gives rise to the {\it monodromy
operator}
$$M_x : H^\cdot (F_x,\Bbb Q) \rightarrow H^\cdot (F_x,\Bbb Q).$$
For any natural number $n$, we consider the {\it Lefschetz number}
$$\Lambda(M^n_x) := \sum_{q \geq 0} (-1)^q \, {\rm Trace} \, (M^n_x,H^q(F_x,\Bbb
Q)),$$
of the $n$-th iterate of $M_x$.  These numbers are related to the {\it monodromy
zeta function of $f$ at $x$}
$$Z^{{\rm mon}}_x (T) := \prod_{q \geq 0}
\, {\rm Det} \, (Id - TM_x,H^q(F_x,\Bbb
Q))^{(-1)^q}$$
as follows : if one writes $\Lambda(M^n_x) = \sum_{i | n} s_i$ for $n \geq 1$,
then $Z^{{\rm mon}}_x (T) = \prod_{i \geq 1} (1  - t^i)^{s_i/i}$.  The
monodromy zeta function of $f$ at $x$ (or equivalently the Lefschetz numbers) is
an important topological invariant of $f$ which has been studied intensively.

\subsubsection{Hodge structures.}\label{3.1.2}
A {\it Hodge structure} is a finite dimensional
$\Bbb Q$-vector\-spa\-ce $H$ together with a bigrading $H \otimes \Bbb C =
\oplus_{p,q \in \Bbb Z} H^{p,q}$, such that $H^{q,p}$ is the complex conjugate
of $H^{p,q}$  and each {\it weight $m$ summand}, $\oplus_{p+q=m} H^{p,q}$, is
defined
over $\Bbb Q$.  The Hodge structures, with the evident notion of morphism, form
an abelian category $\HS$ with tensor product.  The elements of the Grothendieck
group $K_0 (\HS)$ of this abelian category are representable as a formal
difference of Hodge structures $[H] - [H^\prime]$, and  $[H] = [H^\prime]$ iff 
$H \cong H^\prime$.  Note that $K_0 (\HS)$ becomes a ring with respect to the
tensor product.

A {\it mixed Hodge structure} is a finite dimensional $\Bbb Q$-vector space $V$
with
a finite increasing filtration $W_\bullet V$, called the {\it weight
filtration}, such
that the associated graded vector space ${\rm Gr}^W_\bullet(V)$ underlies a Hodge
structure having ${\rm Gr}^W_m(V)$ as weight $m$ summand.  Note that $V$ determines
in a natural way an element [V] in $K_0 (\HS)$, namely $[V] := \sum_m
[{\rm Gr}^W_m(V)]$.

When $X$ is an algebraic variety over $k = \Bbb C$, the simplicial cohomology
groups $H^i_c(X, \Bbb Q)$ of $X$, with compact support, underly a natural mixed
Hodge structure, and the {\it Hodge characteristic} $\chi_h(X)$ of $X$ (with
compact support) is defined by
$$\chi_h(X) := \sum_i (-1)^i [H^i_c(X,\Bbb Q)] \in K_0(\HS).$$
This yields a map $\chi_h :
\Var_\Bbb C \rightarrow K_0 (\HS)$, which is a
generalized Euler characteristic, and which factors through $\mathcal{M}_k$,
because $\chi_k(\Bbb A^1_{k})$ is actually invertible in the ring
$K_0 (\HS)$.  When
$X$ is proper an smooth, the mixed Hodge structure on $H^i_c(X,\Bbb Q)$ is in
fact a Hodge structure, the weight filtration being concentrated in weight $i$.
We refer to \cite{St3} for an introduction to Hodge structures.

\subsubsection{The Hodge spectrum.}\label{3.1.3}
The cohomology groups $H^i_c(F_x,\Bbb Q)$ of
the Milnor fiber $F_x$ carry a natural mixed Hodge structure 
(\cite{St1}, \cite{Sa1}, \cite{Sa2}),
which
is compatible with the monodromy operator $M_x$.  Hence we can define the Hodge
characteristic $\chi_h(F_x)$ of $F_x$ by
$$\chi_h(F_x) := \sum_i (-1)^i [H^i(F_x,\Bbb Q)] \in K_0(\HS).$$
Actually by taking into account the monodromy action we can consider
$\chi_h(F_x)$ as an element of the Grothendieck group $K_0(\HS^{\rm mon})$ of
the abelian category $\HS^{\rm mon}$ of Hodge structures with a quasi-unipotent
endomorphism.  (Quasi-unipotent means that some power of it is unipotent.)
Again $K_0(\HS^{\rm mon})$ is a ring by the tensor product.

There is a natural linear map, called the Hodge spectrum
$$\hsp :
K_0(\HS^{\rm mon }) \rightarrow \Bbb Z[t^{1/\Bbb N}] := \cup_{n
\geq 1} \Bbb
Z[t^{1/n},t^{-1/n}],$$
with $\hsp ([H]) := \sum_{\alpha \in \Bbb Q \cap[0,1[} t^\alpha (\sum_{p,q
\in
\Bbb Z} \dim(H^{p,q})_\alpha)t^p$, for any Hodge structure $H$ with a
quasi-unipotent endomorphism, where $H^{p,q}_\alpha$ is the generalized
eigenspace of $H^{p,q}$ with respect to the eigenvalue
$e^{2\pi\sqrt{-1}\alpha}$.  Note that $\hsp$
is not a ring homomorphism,
although it becomes one when we endow $K_0(\HS^{\rm mon})$ with a different
ring
multiplication, namely the one induced by the operation $\ast$ in section 5.

We recall that $\hsp (f,x) := (-1)^{d - 1}\hsp (\chi_h(F_x)-1)$ is called the Hodge
spectrum of
$f$ at $x$.  It is a very important invariant, with remarkable properties, see
\cite{St1}, \cite{St2}, \cite{Va}.

\subsection{The motivic zeta function}
Let $n \geq 1$ be an integer.  The morphism $f : X \rightarrow \Bbb A^1_k$
induces a morphism $f_n : {\cal L}_n(X) \rightarrow {\cal L}_n (\Bbb A^1_k)$.

Any point $\alpha$ of ${\cal L}(\Bbb A^1_k)$, resp. ${\cal L}_n(\Bbb
A^1_k)$, yields a $K$-rational point, for some field $K$ containing $k$, and
hence a power series $\alpha(t) \in K[[t]],$ resp. $\alpha(t) \in
K[[t]]/t^{n+1}$.    This yields maps
$${\rm ord}_t : {\cal L}(\Bbb A^1_k) \rightarrow \Bbb N \cup \{ \infty \},
\quad {\rm ord}_t : {\cal L}_n(\Bbb A^1_k) \rightarrow \{ 0,1,\cdots,n,\infty
\},$$
with ${\rm ord}_t \alpha$ the largest $e$ such that $t^e$ divides $\alpha(t)$.

We set
$$\frak X_n := \{ \varphi \in {\cal L}_n(X) \, | \, {\rm ord}_t f_n(\varphi) = n
\}.$$
This is a locally closed subvariety of ${\cal L}_n(X)$.  Note that $\frak X_n$
is actually an $X_0$-variety, through the morphism $\pi^n_0 : {\cal L}_n(X)
\rightarrow X$.  Indeed $\pi^n_0(\frak X_n) \subset X_0$, since $n \geq 1$. 
We consider the morphism
$$\bar f_n : \frak X_n \rightarrow \Bbb G_{m,k} := \Bbb A^1_k \setminus \{
0 \},$$
sending a point $\varphi$ in $\frak X_n$ to the coefficient of $t^n$ in
$f_n(\varphi)$.  There is a natural action of $\Bbb G_{m,k}$ on $\frak X_n$
given by $a \cdot \varphi(t) = \varphi(at)$, where $\varphi(t)$ is the vector
of
power series corresponding to $\varphi$ (in some local coordinate system).
Since
$\bar f_n(a \cdot \varphi) = a^n \bar f_n(\varphi)$ it follows
that $\bar f_n$ is a locally trivial fibration.

We denote by $\frak X_{n,1}$ the fiber $\bar f^{-1}_n(1)$.  Note that the
action of $\Bbb G_{m,k}$ on $\frak X_n$ induces a good action of $\mu_n$ (and
hence of $\hat \mu$) on $\frak X_{n,1}$.  Since $\bar f_n$ is a locally
trivial fibration, the $X_0$-variety $\frak X_{n,1}$ and the action of $\mu_n$
on it, completely determines both the variety $\frak X_n$ and the morphism
$$(\bar f_n, \pi^n_0) : \frak X_n \rightarrow  \Bbb
G_{m,k} \times X_0.$$  Indeed it is easy to verify that $\frak X_n$, as a
$(\Bbb
G_{m,k} \times X_0)$-variety, is isomorphic to the quotient of $\frak X_{n,1}
\times \Bbb G_{m,k}$ under the $\mu_n$-action defined by $a(\varphi,b) = (a
\varphi,a^{-1}b)$.

\begin{definition}\label{3.2.1}
The motivic zeta function of $f : X \rightarrow
\Bbb A^1_k$, is the power series over $\mathcal{M}^{\hat \mu}_{X_0}$ defined by
$$Z(T) := \sum_{n \geq 1} \, [\frak X_{n,1}/X_0, \hat \mu] \,
\Bbb L^{-nd} \, T^n.$$
Moreover we define the naive motivic zeta function of $f$ as the power series
over $\mathcal{M}_{X_0}$ defined by
$$Z^{{\rm naive}}(T) := \sum_{n \geq 1} \, [\frak X_n/X_0]
\, \Bbb L^{-nd} \, T^n.$$
\end{definition}

\noindent
Theorem \ref{3.3.1} and Corollary 3.3.2 below show  that $Z(T)$ and
$Z^{{\rm naive}}(T)$
are  rational.  In 3.4 and 3.5  we will see that $Z(T)$ and
$Z^{{\rm naive}}(T)$ give
rise to interesting new
invariants of $f$.  The definition of $Z(T)$ goes back to \cite{DeLo7} (in the
non-relative version working in $\mathcal{M}^{\hat \mu}_k$).  Many related
motivic zeta functions, and $Z^{{\rm naive}}(T)$, were first introduced in
\cite{DeLo2}, inspired by work of Kontsevich \cite{Ko}.  The idea of rather working in
the relative Grothendieck group was introduced by Looijenga \cite{Loo}. 

\subsection{A formula for the motivic zeta function}
We recall that $X_0$ denotes the locus of $f = 0$ in $X$.   Let $(Y,h)$ be a
resolution of $f$.  By this, we mean that $Y$ is a
nonsingular and irreducible algebraic variety over $k, h : Y \rightarrow X$ is
a
proper morphism, that the restriction $h : Y \setminus h^{-1}(X_0) \rightarrow
X
\setminus X_0$ is an isomorphism, and that $h^{-1}(X_0)$ has only normal
crossings
as a subvariety of $Y$. 

We  denote by $E_i, i \in J$, the irreducible  components (over $k$) of
$h^{-1}(X_0)$.  For each $i \in J$, denote by $N_i$ the multiplicity of $E_i$
in
the divisor of $f \circ h$ on $Y$, and by $\nu_i - 1$ the multiplicity of $E_i$
in the divisor of $h^\ast dx$, where $dx$ is a local non vanishing volume form
at any point of $h(E_i)$, i.e. a local generator of the sheaf of differential
forms of maximal
degree.  For $i \in J$ and $I \subset J$, we consider the nonsingular varieties
$E^\circ_i := E_i \setminus \cup_{j \ne i} E_j, E_I = \cap_{i \in I} E_i$, and
$E^\circ_I := E_I \setminus \cup_{j \in J \setminus I} E_j$.

Let $m_I = {\gcd}(N_i)_{i \in I}$.  We introduce an unramified
Galois cover $\tilde E^\circ_I$ of $E^\circ_I$, with Galois group $\mu_{m_I}$,
as follows.  Let $U$ be an affine Zariski open subset of $Y$, such that, on
$U$, $f \circ h = uv^{m_I}$, with $u$ a unit on $U$ and $v$ a morphism from $U$
to $\Bbb A^1_k$.  Then the restriction of $\tilde E^\circ_I$ above $E^\circ_I
\cap U$, denoted by $\tilde E_I^\circ \cap U$, is defined as
$$\{ (z,y) \in \Bbb A^1_k \times (E^\circ_I \cap U) | z^{m_I} = u^{-1} \}.$$
Note that $E^\circ_I$ can be covered by such affine open subsets $U$ of $Y$.
Gluing together the covers $\tilde E^\circ_I \cap U$, in the obvious way, we
obtain the cover $\tilde E^\circ_I$ of $E^\circ_I$ which has a natural
$\mu_{m_I}$-action (obtained by multiplying the z-coordinate with the elements
of $\mu_{m_I}$).  This $\mu_{m_I}$-action on $\tilde E^\circ_I$ induces an
$\hat \mu$-action on $\tilde E^\circ_I$ in the obvious way.

\begin{theorem}[\cite{DeLo7}, \cite{Loo}]\label{3.3.1}
With the previous notations, the
following relation
holds in $\mathcal{M}^{\hat \mu}_{X_0} [[T]]$ :
$$Z(T) =  \sum_{\emptyset \ne I \subset J} \, (\Bbb L - 1)^{|I|-1} \,
[ \tilde
E^\circ_I/X_0, \hat \mu] \,
\prod_{i
\in I} \frac{\Bbb
L^{-\nu_i} T^{N_i}}{1 - \Bbb L^{-\nu_i} T^{N_i}}.$$
\end{theorem}

\noindent

From the above theorem one easily deduces (using e.g. Lemma 5.1 in 
\cite{Loo}) the
following corollary, which is basically a special case of Theorem 2.2.1 in
\cite{DeLo2}.

\begin{corollary}\label{3.3.2}
With the previous notations, the following
relation
holds in \break $\mathcal{M}_{X_0}[[T]]$ :
$$Z^{{\rm naive}}(T) = \sum_{\emptyset \ne I \subset J} \,
(\Bbb L - 1)^{|I|} \,
[E^\circ_I/X_0] \, \prod_{i \in I}
\frac{\Bbb L^{- \nu_i}T^{N_i}}{1 - \Bbb L^{-\nu_i}
\,T^{N_i}}.$$ 
\end{corollary}

\setcounter{subsubsection}{2}

\subsubsection{Proof of the rationality of $J(T)$.}
We defined $J(T)$ in Theorem \ref{2.3.1} above.  We will now discuss the proof of the
rationality of $J(T)$ in the special case when the variety $X$ in Theorem 
\ref{2.3.1}
is the locus $X_0$ of a polynomial $f$ in the affine space $\Bbb A^d_k$.
Let $Z^{{\rm naive}}(T)$ be the naive motivic zeta function of $f : \Bbb A^d_k
\rightarrow \Bbb A^1_k$.  It is straightforward to verify that
$J(T {\Bbb L}^{-d}) = \frac{[X_0]
- Z^{{\rm naive}}(T)}{1-T}$.  Hence the rationality of $J(T)$ is a direct
consequence of Corollary \ref{3.3.2}.

\subsection{The topological zeta functions}
Let $\mathcal{M}_{S,{\rm loc}}$ resp. $\mathcal{M}^{\hat \mu}_{S,{\rm loc}}$,
be the
ring obtained from $\mathcal{M}_S$, resp. $\mathcal{M}^{\hat \mu}_S$, by
inverting the
elements $[\Bbb P^i_k] = 1 + \Bbb L + \Bbb L^2 + \cdots + \Bbb L^i$, for $i =
1,2,3,\cdots$, where $\Bbb P^i_k$ denotes the $i$-dimensional projective space
over $k$.

We keep the notations of 3.3, but take $k = \Bbb C$.  For any integer $s \geq
1$, evaluating $Z^{{\rm naive}}(T)$ at $T = \Bbb L^{-s}$ yields a
well-defined
element of $\mathcal{M}_{X_0,{\rm loc}}$, namely $\sum_{\emptyset \ne I
\subset J} 
[E^\circ_I/X_0] \prod_{i
\in I} [ \Bbb P^{sN_i + \nu_i - 1}]^{-1}$.  Applying the topological Euler
characteristic $\chi_{\top}$ we obtain
\begin{equation}\label{ast}Z_{\top} (s) := \chi_{\top} (Z^{{\rm naive}} (\Bbb L^{-s})) :=
 \underset{\emptyset \ne I \subset J}  \sum  \chi_{\top}
 (E^\circ_I)
\prod_{i \in I} \frac{1}{sN_i + \nu_i}.\tag{$*$}
\end{equation}
We call $Z_{\top}(s)$, considered as a rational function in the variable
$s$, the {\it untwisted topological zeta function} of $f : X \rightarrow \Bbb
A^1_k$.

Evaluating $(\Bbb L - 1)Z(T)$, instead of $Z^{{\rm naive}}(T)$, at $T = \Bbb
L^{-s}$, and applying the equivariant topological Euler characteristic
$\chi_{\top}(-,\alpha)$, with $\alpha : \hat \mu \rightarrow \Bbb C$ a
character of order $e$, we obtain the {\it twisted topological zeta function}
(for any integer $e \geq 1$) :
\begin{equation}\label{2ast}
\begin{split}
Z^{(e)}_{\top}(s)  := & \chi_{\top}((\Bbb L-1)Z(\Bbb
L^{-s}),\alpha) \\
 := & \sum_{\emptyset \ne I \subset J, \, e | m_I} \chi_{\top}(E^\circ_I)
\prod_{i \in I} \frac{1}{sN_i + \nu_i}
\end{split}\tag{$**$}
\end{equation}

Note that, if we would define the topological zeta functions by the
right-hand-side of (\ref{ast}) and (\ref{2ast}), then it would be not at all clear that this is
independent of the choosen resolution.  It is the intrinsic definition using
the motivic zeta function (which is based on the notion of arc spaces) that
makes this independence obvious.  The topological zeta functions were first
introduced by Denef and Loeser in \cite{DeLo1} using $p$-adic integration and the
Grothendieck-Lefschetz trace formula to prove their independence of the
choosen resolution.  Our approach using arc spaces first appeared in 
\cite{DeLo2}.

The topological zeta functions are quite subtle invariants of $f$, and have
been further investigated by Veys \cite{Ve3}, \cite{Ve4}.  There are some fascinating
conjectures about them.

\begin{conjecture}[Monodromy conjecture for $Z^{(e)}_{\top}$]
If $s$ is a pole of $Z^{(e)}_{\top}(s)$ then $e^{2\pi\sqrt{-1}s}$ is
an
eigenvalue of the monodromy action on the cohomology of the Milnor fiber at
some point of the locus of $f$.
\end{conjecture}

\begin{conjecture}[Holomorphy conjecture for $Z^{(e)}_{\top}$]
The function $Z^{(e)}_{\top}(s)$ is a po\-lynomial in $s$, unless
there is
an eigenvalue with order divisible by $e$, of the monodromy action on the
cohomology of the Milnor fiber at some point of the locus of $f$.
\end{conjecture}

\noindent
Loeser \cite{Loe} and Veys \cite{Ve2} proved that these conjectures are true when $X =
\Bbb
A^2_{\Bbb C}$.  A lot of experimental evidence has been obtained by Veys 
\cite{Ve1}
when $X =
\Bbb A^3_{\Bbb C}$.  We refer to \cite{Ve4} and \cite{Ve5} for very interesting
generalizations.

\subsection{Relations with monodromy and the motivic Milnor fiber}
The Lefschetz numbers $\Lambda(M^n_x)$ of $f$ at $x$, which we recalled in
3.1.1 can be expressed in terms of a resolution of $f$, by the following
formula of A'Campo.

\begin{theorem}[A'Campo, \cite{Ac}]\label{3.5.1}
Let $k = \Bbb C$.  Assume the
notations
of 3.1
and 3.3.  Then for any integer $n \geq 0$ we have
$$\Lambda(M^n_x) = \sum_{N_i|n} N_i \, \chi_{\top} (E^\circ_i \cap
h^{-1}(x)).$$
\end{theorem}

\noindent
In particular we see that the right-hand-side of the above formula is
independent of the choosen resolution $h$.  Note that the material in 3.3 and
3.4 yields many other expressions which are independent of the chosen
resolution, but A'Campo's result was probably the first in this direction.

Applying the natural map Fiber$_x : \mathcal{M}^{\hat \mu}_{X_0} \rightarrow
\mathcal{M}^{\hat
\mu}_k :
[A/X_0,\hat \mu] \mapsto [A \times_{X_0} \{ x \}, \hat \mu]$, followed by the
equivariant topological Euler characteristic $\chi_{\top} (-,1)$, on the
coefficients of $Z(T)$ and using Theorem \ref{3.3.1}, we obtain the following
theorem.

\begin{theorem}[\cite{DeLo7}]\label{3.5.2}
Let $k = \Bbb C$, then for any integer $n
\geq
1$ we
have $\Lambda(M^n_x) = \chi_{\top} (\frak X_{n,1,x})$, where
$$\frak X_{n,1,x} := \frak X_{n,1} \times_{X_0} \{ x \}.$$
\end{theorem}

\noindent
Thus we see that the monodromy zeta function of $f$ at $x$ is completely
determined by the motivic zeta function $Z(T)$.  Next, we will see that also
the Hodge spectrum of $f$ at $x$ is determined by $Z(T)$.

\begin{definition}[\cite{DeLo2}, \cite{DeLo7}]\label{3.5.3}
Expanding the rational function $Z(T)$ as a power series in $T^{-1}$ and taking
minus its constant term, yields a well defined element of $\mathcal
{M}_{X_0}^{\hat
\mu}$, namely
$$\mathcal{S} := - \lim_{T \rightarrow \infty} Z(T) := \sum_{\emptyset \ne I
\subset J} (1 - \Bbb L)^{|I|-1} [ \tilde E^\circ_I].$$
Moreover we set $\mathcal{S}_x := {\rm Fiber}_x(\mathcal{S}) \in
\mathcal{M}^{\hat
\mu}_k$.  Instead of $\mathcal{S}$ and $\mathcal{S}_x$ we will also write
$\mathcal{S}_f$ and $\mathcal{S}_{f,x}$.  These definitions hold for any field
$k$ of characteristic zero.
\end{definition}

\noindent
Note again that the most right-hand-side of the above formula is independent of
the choosen resolution (because of its relation to $Z(T)$), although a priori
this is not at all evident.

We strongly believe that $\mathcal{S}_x$ is the correct virtual motivic
incarnation of
the Milnor fiber $F_x$ of $f$ at $x$ (which is in itself not at all motivic).
We will see below (Theorem \ref{3.5.5}) that this is indeed true for the Hodge
realization. A similar result holds for $\ell$-adic cohomology, see 
\cite{De3}.
Moreover we strongly believe that $\mathcal{S}$ is the virtual
motivic
incarnation of the so called complex of nearby cycles $\psi_f$ of $f$, which is
a complex of sheaves on $X_0$.  For the definition of $\psi_f$ and the complex
$\phi_f$ of vanishing cycles, we refer to \cite{SGA 7}, Exp. XIII; but we will not need these
notions in the present survey.  Inspired by the notation $\phi_f$ from the
theory of vanishing cycles, we introduce the following

\begin{notation}\label{3.5.4}We set $\mathcal{S}^\phi_f := 
(-1)^{d - 1} (\mathcal{S}_f - [X_0])
\in \mathcal{M}^{\hat \mu}_{X_0}$
and
$\mathcal{S}^\phi_{f,x} := (-1)^{d - 1} (\mathcal{S}_{f, x}
- 1) \in \mathcal{M}^{\hat \mu}_k$.
\end{notation}

\noindent
We regard $\mathcal{S}^\phi_f$ as the virtual motivic incarnation of the
complex $\phi_f [d - 1]$.

Assume now again that $k = \Bbb C$.  We denote by $\chi_h$ the canonical ring
homomorphism (called the Hodge characteristic)
$$\chi_h : \mathcal{M}^{\hat \mu}_k \rightarrow K_0({\HS^{\rm mon}}),$$
which associates to any complex algebraic variety $Z$, with a good
$\mu_n$-action, its Hodge characteristic as defined in 3.1.2, together with the
endomorphism induced by $Z \rightarrow Z : z \mapsto e^{2 \pi\sqrt{-1}/n}z$.
(For the definition of $K_0(\HS^{\rm mon})$, see 3.1.3.)

\begin{theorem}[\cite{DeLo2}]\label{3.5.5}
Assume the above notation
with $k
= \Bbb C$, and the notation of 3.1.  Then we have the following equality in
$K_0({\HS^{\rm mon}})$ :
$$\chi_h(F_x) = \chi_h(\mathcal{S}_x).$$
\end{theorem}

\noindent
Moreover this theorem can be enhanced as an equality in the Grothendieck group
of the abelian category of variations of Hodge structures with a
quasi-unipotent endomorphism, when we replace $\mathcal{S}_x$ by $\mathcal{S}$,
and $F_x$ by $\psi_f$.

\noindent
Theorem \ref{3.5.5} yields that $\hsp(f,x) =
{\rm hsp}(\chi_h(\mathcal{S}^\phi_{f,x}))$.  Thus the motivic zeta function
$Z(T)$ completely determines the Hodge spectrum of $f$ at $x$.

\section{Motivic integration and the proof of Theorem \ref{3.3.1}}\label{sec4}
The notion of motivic integration on $\mathcal L(X)$ is due to Kontsevich 
\cite{Ko},
who discovered its basic properties when $X$ is nonsingular.  This subject has
been  further developed by Batyrev \cite{Ba2}, \cite{Ba3} and Denef-Loeser
\cite{DeLo2}, \cite{DeLo3}, \cite{DeLo4}, \cite{DeLo5}, \cite{DeLo6},
\cite{DeLo7}.
See also the recent report by
Looijenga \cite{Loo} which contains some substantial improvements.  Actually the
best way to understand motivic integration is to consider it as being an 
analogue of
$p$-adic integration, cf. section 6.

Let $X$ be an algebraic variety over $k$ of pure dimension $d$, not necessarily
nonsingular.  Let $X_{{\rm sing}}$ denote the singular locus of $X$.

\subsection{Naive motivic integration}
A subset of $A$ of ${\cal L}(X)$ is called {\it constructible} if $A =
\pi^{-1}_n(C)$ with $C$ a constructible subset of ${\cal L}_n(X)$ for some
integer $n \geq 0$.   A subset $A$ of ${\cal L}(X)$ is called {\it stable} if it
is constructible and $A \cap {\cal L}(X_{{\rm sing}}) = \emptyset$.  If $A
\subset {\cal L}(X)$ is stable, then $[\pi_n(A)] \Bbb L^{-(n+1)d}$, considered
as an
element of $\mathcal{M}_k$, stabilizes for $n$ big enough, and
$$\tilde \mu(A) := \lim_{n \rightarrow \infty} [\pi_n(A)]\Bbb L^{-(n+1)d} \in
\mathcal{M}_k$$
is called the {\it naive motivic measure} of $A$.  When $X$ is nonsingular,
this claim follows from the fact that the natural maps ${\cal L}_{n+1}(X)
\rightarrow {\cal L}_n(X)$ are locally trivial fibrations with fiber
$\Bbb A^d_{k}$.
In the general case, the claim follows from \cite{DeLo3}, Lemma 4.1.

When $\theta : A \rightarrow \mathcal{M}_k$ is a map with finite image whose
fibers are stable subsets of ${\cal L}(X)$, we define the integral $\int_A
\theta d \tilde \mu := \sum_{c \in {\rm Image} \, \theta} c \tilde \mu
(\theta^{-1}(c))$.  The most fundamental result in the theory of arc spaces is
the following change of variables formula, which was first obtained by
Kontsevich \cite{Ko} when $X$ is nonsingular.

\begin{theorem}[\cite{Ko}, \cite{DeLo3}, \cite{DeLo5}]\label{4.1.1}
Let $h : Y \rightarrow X$ be a
morphism of
algebraic varieties over $k$.  Suppose that $h$ is birational and proper.  Let
$A \subset {\cal L}(X)$ be stable and suppose that
${\rm ord}_t{\rm Jac}_h$ is bounded on
$h^{-1}(A) \subset {\cal L}(Y)$. Then
$$\tilde \mu (A) = \int_{h^{-1}(A)} \Bbb L^{- {\rm ord}_t{\rm Jac}_h} d
\tilde
\mu.$$
\end{theorem}

\noindent
In the above theorem, ${\rm ord}_t{\rm Jac}_h$, for $y \in {\cal L}(Y)$, denotes the
$t$-order of the Jacobian of $h$ at $y$.  When $X$ and $Y$ are nonsingular this
is the ord$_t$ of the determinant of the Jacobian matrix of $h$ at $y$ with
respect to any system of local coordinates on $X$ and on $Y$.  For the
definition of ${\rm ord}_t{\rm Jac}_h$, in the general
case, we refer to \cite{DeLo3} and
\cite{DeLo5}.

\subsection{About the proof of Theorem \ref{3.3.1}}
The proof of Theorem \ref{3.3.1} consists of an explicit calculation of $[\frak
X_{n,1}/X_0, \hat \mu] \in \mathcal{M}^{\hat \mu}_{X_0}$ for each $n$.  Note
that in $\mathcal{M}_k$ we have the equality
$$[\frak X_{n,1}] = \Bbb L^{(n+1)d} \tilde \mu (\pi^{-1}_n(\frak X_{n,1})).$$
Thus using the change of variables formula 4.1.1, we see that $[\frak X_{n,1}]$
is equal to an integral over a stable subset of ${\cal L}(Y)$, where $h : Y
\rightarrow X$ is a resolution of $f$ as in 3.3.  Because $f \circ h$ is
locally a monomial, that integral can be explicitely calculated and yields an
explicit expression for $[\frak X_{n,1}]$ as an element of $\mathcal{M}_k$.
Taking
into account the
$\mu_n$-action on $\frak X_{n,1}$ and the natural map $\frak X_{n,1}
\rightarrow X_0$, one actually obtains a similar formula for $[\frak
X_{n,1}/X_0,\hat \mu]$, which yields Theorem \ref{3.3.1}.

\subsection{Motivic integration}
Let $A$ be a constructible subset of ${\cal L}(X)$.  When $A$ is not stable,
$[\pi_n(A)] \Bbb L^{-(n+1)d}$ will not always stabilize.  However it is easy to
prove (see \cite{DeLo3}) that the limit
$$\mu(A) := \lim_{n \rightarrow \infty} [\pi_n(A)] \Bbb L^{-(n+1)d}$$
exists in the {\it completed Grothendieck group} $\hat \mathcal{M}_k$, which is
the completion of $\mathcal{M}_k$ with respect to the filtration $F^m
\mathcal{M}_k, m \in \Bbb Z$, where $F^m \mathcal{M}_k$ is the subgroup of
$\mathcal{M}_k$ generated by the elements $[S] \Bbb L^{-i}$, with $S \in
\Var_k$, $i - \dim S \geq m$.  The completed Grothendieck ring
$\mathcal{M}_k$ was first introduced by Kontsevich.  (In a similar way one can
define the completions $\hat \mathcal{M}_S$ and $\hat \mathcal{M}^{\hat \mu}_S$
of $\mathcal{M}_S$ and $\mathcal{M}^{\hat \mu}_S$.)  The element $\mu(A)$ of
$\hat \mathcal{M}_k$ is called the {\it motivic measure of} $A$.  This yields a
$\sigma$-additive measure $\mu$ on the Boolean algebra of constructible subsets
of ${\cal L}(X)$.  Actually all the above still works when $A$ is a
semi-algebraic subset of ${\cal L}(X)$, cf. \cite{DeLo3}.  It is even possible to
define
the notion of a measurable subset of ${\cal L}(X)$ and to integrate measurable
functions on ${\cal L}(X)$, see \cite{Ba2}, \cite{DeLo5}.

The change of variables formula 4.1.1 remains true with $\tilde \mu$ replaced
by
$\mu$, for any constructible (or measurable) subset of ${\cal L}(X)$, without
assuming that ${\rm ord}_t{\rm Jac}_h$
is bounded on $h^{-1}(A)$.

It is not known whether the natural map $\mathcal{M}_k \rightarrow \hat
\mathcal{M}_k$ is injective, but the topological Euler characteristic, the
Hodge-Deligne polynomial, the Hodge
characteristic,
and many other important generalized Euler
characteristics all factor through the image $\bar \mathcal{M}_k$ of
$\mathcal{M}_k$ in $\hat \mathcal{M}_k$ (after inverting the image of $\Bbb L$
in the target ring).

We can  consider the motivic volume of the whole arc space ${\cal L}(X)$, namely
$\mu({\cal L}(X))$.  Clearly, when $X$ is nonsingular, $\mu({\cal L}(X)) = [X]
\Bbb L^{-d}$ in $\hat \mathcal{M}_k$.  Here and in what follows, we denote the
image of $[X]$, resp. $\Bbb L$, in $\hat \mathcal{M}_k$ again by $[X]$, resp.
$\Bbb L$.  When $X$ is not necessarily nonsingular,  we can calculate
$\mu({\cal L}(X))$ using a
suitable resolution of singularities $h : Y \rightarrow X$ of $X$.  More
precisely we have the following

\begin{theorem}\label{4.3.1}
Let $h : Y \rightarrow X$ be a proper birational
morphism
with $Y$ nonsingular.  Assume that the exceptional locus of $h$ has normal
crossings and that the image of $h^\ast(\Omega^d_X)$ in $\Omega^d_Y$ is an
invertible sheaf, where $\Omega^d_X$ and $\Omega^d_Y$ denote the sheaf of
differential forms of maximal degree.  Let $E_j, j \in J$, be the
$k$-irreducible components of the exceptional locus of $h$.  For any subset $I$
of $J$, set $E^\circ_I = (\cap_{i \in I} E_i) \setminus \cup_{j \in J \setminus
I} E_j$.  For $i \in I$, let $\nu_i - 1$ be the multiplicity along $E_i$ of the
divisor associated to
$h^\ast(\Omega^d_X)$.  Then, in $\hat \mathcal{M}_k$, we have
$$\mu({\cal L}(X)) = \Bbb L^{-d} \sum_{I \subset J}
\,
[E^\circ_I]
\,
\prod_{i \in I} \,
[\Bbb P^{\nu_i - 1}]^{-1}.$$
\end{theorem}

\noindent
In particular we see that $\mu({\cal L}(X)) \in \bar \mathcal{M}_{k,{\rm loc}}
\subset \hat \mathcal{M}_k$, where $\bar \mathcal{M}_{k,{\rm loc}}$ denotes
the
ring obtained from $\bar \mathcal{M}_k$ by inverting the elements $1 + \Bbb L +
\cdots + \Bbb L^i$, for all $i = 1,2,3,\cdots$.

\noindent
About the proof of this theorem, we remark that
$\mu({\cal L}(X)) = \int_{{\cal L}(Y)} \Bbb L^{- {\rm ord}_t{\rm Jac}_h}d \mu$, by the change of variables
formula.  Because ${\rm Jac}_h$ is locally a monomial, this integral can be easily
calculated, which yields the theorem.

\subsection{Applications}

\subsubsection{New invariants of singular varieties.}
Suppose $k = \Bbb C$.  Since $\chi_{\top}$ and $\chi_{hp}$ factor through
$\bar \mathcal{M}_k$, we have natural maps $\chi_{\top} : \bar
\mathcal{M}_{k, {\rm loc}} \rightarrow \Bbb Q$ and $\chi_{hp} : \bar
\mathcal{M}_{k,{\rm loc}} \rightarrow \Bbb Z[[u,v]][u^{-1},v^{-1}]$.  Hence we
can consider $\chi_{\top}({\cal L}(X)) \in \Bbb Q$ and $(uv)^d
\chi_{{\rm hp}}({\cal L}(X)) \in \Bbb Z[[u,v]]$, which are new invariants of
$X$ when $X$ is singular.  When $X$ is nonsingular, these invariants equal
$\chi_{\top}(X)$, resp. $\chi_{{\rm hp}}(X)$.  We call the coefficients
of $(uv)^d \chi_{{\rm hp}}({\cal L}(X))$ (with an appropriate sign change) the
arc-Hodge numbers of $X$.  When $X$ has only canonical Gorenstein
singularities, Batyrev \cite{Ba2} introduced the so called {\it stringy Hodge
numbers}
of $X$, which are obtained in a similar way, replacing $\mu({\cal L}(X))$ by
$\int_{{\cal L}(X)} \Bbb L^{- {\rm ord}_t \omega_X} d \mu$, where $\omega_X$ denotes
the canonical class of $X$.  The stringy Hodge numbers play an important role
in the work of Batyrev on mirror symmetry, see \cite{Ba2}, \cite{Ba4}, 
\cite{BaDi}, \cite{BaBo}.
Other fascinating related invariants were obtained by Veys \cite{Ve4}, 
\cite{Ve5}.
\subsubsection{Calabi-Yau manifolds.}
Let $X$ and $Y$ be two Calabi-Yau manifolds, i.e. nonsingular proper complex
algebraic varieties which admit a nonvanishing differential form of maximal
degree, which we denote respectively by $\omega_X$ and $\omega_Y$.  Kontsevich
\cite{Ko} proved that $X$ and $Y$ have the same Hodge numbers and the same Hodge
structure on their cohomology, when $X$ and $Y$ are birationally equivalent.
The proof goes as follows : There exists a nonsingular proper complex algebraic
variety $Z$ and birational morphisms $h_X : Z \rightarrow X$ and $h_Y : Z
\rightarrow Y$.  Note that $(h_Y \circ h^{-1}_X)^\ast (\omega_Y)$ equals $c \
\omega_X$  for
some $c \in  \Bbb C^\times$ because $\omega_X$ has no zeroes.  Hence
$c \ h^\ast_X(\omega_X) = h^\ast_Y(\omega_Y)$.  Thus
${\rm ord}_t{\rm Jac}_{h_X}
=
{\rm ord}_t{\rm Jac}_{h_Y}$
on ${\cal L}(Z)$, and by the change of variables formula both
$\mu({\cal L}(X))$ and $\mu({\cal L}(Y))$ equal the same integral on
${\cal L} (Z)$.  Because $\mu({\cal L}(X)) = [X] \Bbb L^{-d}$ and $\mu({\cal L}(Y)) =
[Y] \Bbb L^{-d}$, this implies that $[X] = [Y]$ in $\bar \mathcal M_k$, which
finishes the proof.

Actually Batyrev \cite{Ba1} first proved that 
$X$ and $Y$ have the same Betti numbers
using $p$-adic integration and the Weil conjectures, and Kontsevich
invented motivic
integration to prove that $X$ and $Y$ have the same Hodge numbers.

\subsubsection{Euler characteristics and modifications.}
Let $h : Y \rightarrow X$ be a modification of nonsingular algebraic varieties
over $k$, meaning that $h$ is a proper birational morphism.  Assume that the
exceptional locus of $h$ has normal crossings, and let $J, E_i, E^\circ_I$ and
$\nu_i$ be as in Theorem \ref{4.3.1}.  Because $X$ is nonsingular,
$\mu({\cal L}(X)) =
[X] \Bbb L^{-d}$ and Theorem \ref{4.3.1} yields the following equality in $\bar
\mathcal{M}_{k,{\rm loc}}$ :
\begin{equation}\label{anotherstar}
[X] = \sum_{I \subset J} \, [E^\circ_I] \, \prod_{i \in I} \,[ \Bbb
P^{\nu_i-1}]^{-1}. \tag{$*$}
\end{equation}

a) When $k = \Bbb C$, applying the topological Euler characteristic on 
(\ref{anotherstar})
yields $\chi_{\top}(X) = \sum_{I \subset J} \chi (E^\circ_I)/\prod_{i \in
I} \nu_i$.  This surprising formula about the Euler characteristic of
modifications was first obtained in \cite{DeLo1} using $p$-adic integration and the
Grothendieck-Lefschetz trace formula.

b) When $k = \Bbb Q$, applying the conductor (with respect to $\ell$-adic cohomology,
see section 2.2) yields the following remarkable formula for the conductor
$c(X)$ of $X$ :
$$c(X) = \prod_{I \subset J} c(E^\circ_I)^{1/\prod_{i \in I} \nu_i}.$$

\section{The motivic Thom-Sebastiani Theorem}\label{sec5}
Let $k$ be a field of characteristic zero, $X$ and $Y$ nonsingular irreducible
algebraic varieties over $k$, $f: X \rightarrow \Bbb A^1_k, \ g: Y \rightarrow
\Bbb A^1_k$ non constant morphisms, and $x \in X(k), y \in Y(k)$.  We denote by
$f \ast g$ the morphism
$$f \ast g : X \times Y \rightarrow \Bbb A^1_k : (x,y) \mapsto f(x) + g(y).$$
The following theorem was first proved by A. Varchenko \cite{Va}, when $x$ and $y$ are
isolated singular points of $f$ and $g$, and by M. Saito \cite{Sa3}, 
\cite{Sa4} in the
general
case.  A similar but much weaker result for the eigenvalues of monodromy was
first proved by Thom and Sebastiani \cite{SeTh}.
\smallskip

\noindent
{\bf Theorem 5.1} (Thom-Sebastiani Theorem for the Hodge 
spectrum){\bf .} {\it
Assume
the notation of 3.1.3, with $k = \Bbb C$.  We have the following equality in
$\Bbb Z[t^{1/\Bbb N}]$ :}
$$\hsp(f \ast g,(x,y)) = \hsp(f,x) \, \hsp(g,y).$$

\noindent
We recall that $\hsp(f,x) = \hsp(\chi_h(\mathcal{S}^\phi_{f,x}))$, with
the
notation of 3.5.4.  We will see next that the above theorem is a direct
consequence of a much stronger result which expresses $\mathcal{S}^\phi_{f
\ast g,(x,y)}$ in terms of $\mathcal{S}^\phi_{f,x}$ and
$\mathcal{S}^\phi_{g,y}$.  Below, we define a binary operation $\ast$ on
$\mathcal{M}^{\hat \mu}_k$ which yields an alternative ring
structure on $\mathcal{M}^{\hat \mu}_k$, such that $\hsp \circ \chi_h$
becomes a
homomorphism of rings (which is not true for the usual multiplication on
$\mathcal M^{\hat \mu}_k$).

Using the theory of arc spaces and the definition of $\mathcal{S}_f$ in terms
of the motivic zeta function $Z(T)$,
work of Denef, Loeser and Looijenga yields the
following theorem
\smallskip

\noindent
{\bf Theorem 5.2} (Motivic Thom-Sebastiani Theorem){\bf .} {\it
Let
$k$ be
a
field of characteristic zero.  Then $\mathcal{S}^\phi_{f \ast g,(x,y)}$
and $\mathcal{S}^\phi_{f,x} \ast \mathcal{S}^\phi_{g,y}$ are equal in
$\mathcal{M}^{\hat \mu}_k$, where the operation $\ast$ is
defined below.}

\smallskip
\noindent
Actually Denef and Loeser \cite{DeLo4}
first proved the above equality in the completed
Gro\-then\-dieck group of Chow motives.  Later Looijenga \cite{Loo}
introduced the
operation
$\ast$ and proved, using basically the same method,
an equality which is similar to Theorem 5.2.
The proof of Theorem 5.2 uses arc spaces in a very essential way
by deriving first a formula relating the motivic zeta functions of $f
\ast g$, $f$ and $g$, and taking afterwards the limit for $T
\rightarrow \infty$.
More precisely, set
$$
Z^{\phi}_{f} (T) :=
(-1)^{d- 1}
\Bigl[Z_f (T) + [X_0] + \frac{Z^{{\rm naive}}_f (T) - [X_0]}{1 - T}\Bigr],
$$
where $Z^{\phi}_{f} (T)$, resp. $Z^{{\rm naive}}_f (T)$, is the
motivic, resp. naive motivic, zeta function
of $f$, $X_0$ is the locus of $f = 0$ in $X$, and $d$ is the dimension
of $X$.
Let $Z^{\phi}_{f, x} (T)$ be obtained from $Z^{\phi}_{f} (T)$
by applying the map ${\rm Fiber}_x$ to its coefficients.
Clearly 
$- \lim_{T \rightarrow \infty} Z^{\phi}_{f} (T)$
is $\mathcal{S}_f^{\phi}$. One proves that
$$
Z^{\phi}_{f \ast g, (x, y)} (T) =
Z^{\phi}_{f, x} (T) \ast
Z^{\phi}_{g, y} (T)
$$
in 
$\mathcal{M}^{\hat \mu}_k$, where $\ast$ is defined coefficientswise.
This implies Theorem 5.2 because
$- \lim_{T \rightarrow \infty}$ commutes with $\ast$ on such power
series without constant terms.

Finally, we explain the definition of the operation $\ast$ on
$\mathcal{M}^{\hat \mu}_k$.  Let $X$ and $Y$ be algebraic varieties over $k$
with good $\mu_n$-action, for some integer $n \geq 1$.  Let $J_n$ be the Fermat
curve in $(\Bbb A^1_k \setminus \{ 0 \})^2$ defined by $u^n + v^n = 1$.  There
is an action of $\mu_n \times \mu_n$ on $J_n$ given by $(\xi,\xi^\prime) \cdot
(u,v) := (\xi u, \xi^\prime v)$.    We define $J_n(X,Y)$
in $\Var_k$ as the quotient of $J_n \times X \times Y$ under the equivalence
relation given by $(\xi u, \xi^\prime v,x,y) = (u,v, \xi x,\xi^\prime y)$ for
all $\xi,\xi^\prime \in \mu_n$.  We let $\mu_n$ act on $J_n(X,Y)$ by $\xi \cdot
(u,v,x,y) := (\xi u, \xi v, x, y)$.  This yields an element $[J_n(X,Y)]$ in
$\mathcal{M}^{\hat \mu}_k$.  If $m$ is a divisor of $n$, and the action of
$\mu_n$ on $X$ and $Y$ factors through $\mu_m$, then $J_m(X,Y) = J_n(X,Y)$.
Thus, in this way, we obtain a binary operation $J : \mathcal{M}^{\hat \mu}_k
\times \mathcal{M}^{\hat \mu}_k \rightarrow \mathcal{M}^{\hat \mu}_k$, which
was first introduced by Looijenga \cite{Loo}.  The operation $J$ is commutative and
bilinear over $\mathcal{M}_k$, considering $\mathcal{M}^{\hat \mu}_k$ as a
module over $\mathcal{M}_k$ through the natural map $\mathcal{M}_k \rightarrow
\mathcal{M}^{\hat \mu}_k$.  One verifies that $J(a,1) = (\Bbb L - 1)\bar a -
a$, where $a \mapsto \bar a : \mathcal{M}^{\hat \mu}_k \rightarrow
\mathcal{M}_k$ is the morphism induced by $[Z] \rightarrow$ [space of $\hat
\mu$-orbits of $Z]$, for any $k$-variety $Z$ with good $\hat 
\mu$-action, cf. \cite{Loo}. In particular we see that 1 is not a neutral element for the operation
$J$.  For this reason it is natural to introduce the operation $\ast$ on
$\mathcal{M}^{\hat \mu}_k$ given by    
$$a \ast b = -J(a,b) + (\Bbb L - 1)\overline{ab},$$
for
$a$ and $b$ in
$\mathcal{M}^{\hat \mu}_k$.
Clearly the operation $\ast$ is commutative and bilinear over $\mathcal{M}_k$,
and $a \ast 1 = a$ for all  $a$ in $\mathcal{M}^{\hat \mu}_k$.  Moreover one
easily verifies that
$\hsp \circ \chi_h$ is a ring
homomorphism with respect to the alternative ring structure
on $\mathcal{M}^{\hat \mu}_k$ given by $\ast$.

\section{The arithmetic motivic Poincar{\'e} series $P_{{\rm arith}}(T)$}\label{sec6}
\subsection{The $p$-adic case}
Assume that $X$ is an algebraic variety over $\Bbb Z$, i.e. a reduced separated
scheme of finite type over $\Bbb Z$.  Let $p$ be a prime number.  We consider
the Poincar{\'e} series
$$J_p(T) = \sum_{n \in \Bbb N} \# X(\Bbb Z/p^{n +1} \Bbb Z)\,T^n, \quad P_p(T) =
\sum_{n
\in \Bbb N} \# (\pi_n(X(\Bbb Z_p)))\,T^n,$$
where $\Bbb Z_p$ denotes the ring of $p$-adic integers and $\pi_n$ is the
natural projection $\pi_n : \Bbb Z_p \rightarrow \Bbb Z/p^{n +1} \Bbb Z$.  Igusa
\cite{Ig}, resp. Denef \cite{De1}, proved that $J_p(T)$, resp. $P_p(T)$, is a rational
function of $T$.  The proofs are based on $p$-adic integration, resolution of
singularities, and for $P_p(T)$ also the theory of $p$-adic semi-algebraic
sets.  Actually the proof of Theorem \ref{2.3.1} about the rationality of $J(T)$ and
$P(T)$ was very much inspired by the proofs
of the rationality of $J_p(T)$ and $P_p(T)$, replacing $p$-adic integration by
motivic integration.  As a matter of fact, for
all the material discussed in the
previous sections, $p$-adic counterparts exist which
were 
discovered first, see \cite{De2} for a survey.   
\subsection{Comparing $J(T)$ and $J_p(T)$}
For any rational power series $G(T)$ over $K_0(\Var_{\Bbb Q})$ (with
denominator
a product of polynomials of the form $1 - \Bbb L^a T^b, b \in \Bbb N \setminus
\{ 0 \}, a \in \Bbb Z$) we choose representatives in $K_0(\Var_{\Bbb Z})$
for
the coefficients in $K_0(\Var_k)$ of numerator and denominator.  In this
way we
find a power series over $K_0(\Var_{\Bbb Z})$, and, for any prime number
$p$, we
can apply to each coefficient the operation $N_p : K_0(\Var_{\Bbb Z})
\rightarrow \Bbb Z : [X] \mapsto \# X(\Bbb Z/p \Bbb Z)$.  This yields a power
series over $\Bbb Z$ which we will denote by $N_p(G(T))$.  If we choose other
representatives in $K_0(\Var_{\Bbb Z})$, the resulting power series
$N_p(G(T))$
will be the same for almost all $p$ (i.e. for all but finitely many prime
numbers $p$).

Comparing the proof of the rationality of $J(T)$ and $J_p(T)$ actually yields
the following

\begin{theorem}\label{6.2.1}
Assume the notation of 6.1 and 6.2.  For almost all
$p$
we have $J_p(T) = N_p(J(T))$.
\end{theorem}

\noindent
Also the motivic zeta functions $Z(T)$ and $Z^{\rm naive}(T)$, have similar
arithmetic interpretations, related to Igusa's local zeta functions, see
\cite{DeLo2}.
However it is not true in general that $N_p(P(T)) = P_p(T)$ for almost all $p$. 
Indeed $N_p(P(T))$ does not count the elements of $X(\Bbb Z/p^{n +1} \Bbb Z)$ which
can be lifted to $X(\Bbb Z_p)$, but counts (for almost all $p$) the elements of
$X(\Bbb Z/p^{n +1} \Bbb Z)$ which can be lifted to $X(\Bbb Z^{{\rm unram}}_p)$, where
$\Bbb Z^{{\rm unram}}_p$ is the maximal unramified extension of $\Bbb
Z_p$.
Note
that the residue field of $\Bbb Z^{{\rm unram}}_p$ is the algebraic closure of $\Bbb
Z/p \Bbb Z$.
\subsection{The motivic Poincar{\'e} series $P_{{\rm arith}}(T)$}
The above discussion leads to the question of defining in a canonical way a
power series $P_{{\rm arith}}(T)$ over (some localization of) $K_0(\Var_Q)$ such
that
$N_p(P_{{\rm arith}})(T) = P_p(T)$ for almost all $p$.

In our recent paper \cite{DeLo6}, we construct, for any algebraic variety over $k$,
in a
canonical way, a rational power series $P_{{\rm arith}}(T)$ over
$K_0({\rm Mot}_k)
\otimes
\Bbb Q$, such that if $k = \Bbb Q$ then $N_p(P_{{\rm arith}}(T)) = P_p(T)$ for almost all
$p$.  Here ${\rm Mot}_k$ denotes the category of Chow motives over $k$.  We
refer to
\cite{Sc} for the definition of this important category, and we only remark here
that
there is a natural ring morphism $K_0(\Var_k) \rightarrow K_0({\rm
  Mot}_k)$,
see \cite{GS}.
Actually
the coefficients of $P_{{\rm arith}}(T)$ are in the image of $K_0(\Var_k) \otimes
\Bbb
Q$.  We need to work at the level of Chow motives, to make our construction
canonical.

The proof of our result is rather complicated and uses several results from
mathematical logic (quantifier elimination for valued fields and finite
fields). For a survey on such relations between logic, geometry and arithmetic,
we refer to \cite{De4}. Examples seem to suggest that $P_{{\rm arith}}(T)$ captures more
geometric
information than $P(T)$, but very little is presently known about it !


\begin{thebibliography}{DeLo 7}

\bibitem[Ac] {Ac}N. A'Campo,
\textit{La fonction z{\^e}ta d'une monodromie},
Comment. Math. Helv.
\textbf{50}
(1975),
233--248.
\bibitem[BaBo] {BaBo}V. Batyrev, L. Borisov, \textit{Mirror duality and string-theoretic Hodge
numbers}, Invent. Math. \textbf{126} (1996), 183--203.
\bibitem[BaDi] {BaDi}V. Batyrev, D. Dais, \textit{Strong McKay correspondence,
string-theoretic Hodge numbers and mirror symmetry}, Topology
\textbf{35}
(1996), 901--929. 
\bibitem[Ba1] {Ba1}V. Batyrev, \textit{Birational Calabi-Yau n-folds have equal
Betti numbers}, in New trends in algebraic geometry, Klaus Hulek et al., eds.,
CUP (1999), 1--11.
\bibitem[Ba2] {Ba2}V. Batyrev, \textit{Stringy Hodge numbers of varieties with
Gorenstein canonical singularities}, in Integrable systems and algebraic
geometry (Kobe/Kyoto, 1997), 1--32, World Sci. Publishing, River Edge, NJ, 1998.
\bibitem[Ba3] {Ba3}V. Batyrev,
\textit{Non-archimedian integrals and stringy Euler numbers of log
terminal pairs}, Journal of European Math. Soc. \textbf{1} 
(1999), 5--33.
\bibitem[Ba4] {Ba4}V. Batyrev, \textit{Mirror symmetry and toric geometry}, Doc.
Math., J. DMV, Extra Vol. ICM Berlin 1998, {\bf II} (1998), 239--248.
\bibitem[De1] {De1}J. Denef, \textit{On the rationality of the Poincar{\'e} series
associated to the $p$-adic points on a variety}, Invent. Math., \textbf{77}
(1984), 1--23.
\bibitem[De2] {De2}J. Denef, {\it Report on Igusa's local zeta function}, in
S{\'e}minaire Bourbaki, volume 1990/91, expos{\'e} 741.  Ast{\'e}risque {\bf
201-202-203} (1991), 359--386.
\bibitem[De3] {De3}J. Denef, {\it Degree of local zeta functions and monodromy},
Compositio Math., {\bf 89} (1993), 207-216.
\bibitem[De4] {De4}J. Denef, {\it Arithmetic and Geometric Applications of
Quantifier Elimination for Valued Fields}, Introduction to the Model Theory of
Fields : Proceedings of the MSRI Workshop 1998, ed. D. Haskell, A. Pillay, and
C. Steinhorn, Cambridge University Press (to appear, 25 pages).
\bibitem[DeLo1] {DeLo1}J. Denef, F. Loeser, {\it Caract{\'e}ristiques
d'Euler-Poincar{\'e}, fonctions z{\^e}ta locales et modifications
analytiques}, J.
Amer. Math. Soc.,  {\bf 5} (1992), 705-720.
\bibitem[DeLo2] {DeLo2}J. Denef, F. Loeser,
\textit{Motivic Igusa zeta functions},
J. Algebraic Geom.,
\textbf{7}
(1998),
505--537.
\bibitem[DeLo3] {DeLo3}J. Denef, F. Loeser,
\textit{Germs of arcs on singular algebraic varieties
and motivic integration},
Invent. Math.
\textbf{135}
(1999),
201--232.
\bibitem[DeLo4] {DeLo4}J. Denef, F. Loeser,
\textit{Motivic exponential integrals and a motivic Thom-Sebastiani
Theorem},
Duke Math. J., 
\textbf{99}
(1999),
285--309.
\bibitem[DeLo5] {DeLo5}J. Denef, F. Loeser, {\it Motivic integration, quotient
singularities and the McKay correspondence}, 20 p., available at math.
AG/9903187.
\bibitem[DeLo6] {DeLo6}J. Denef, F. Loeser, {\it Definable sets, motives and $P$-adic
integrals}, 45 p., available at math. AG/9910107.
\bibitem[DeLo7] {DeLo7}J. Denef, F. Loeser, {\it Lefschetz numbers of the monodromy
and truncated arcs}, 10 p., available at math. AG/0001105.
\bibitem[GS] {GS}H. Gillet, C. Soul{\'e},
\textit{Descent, motives and $K$-theory},
J. Reine Angew. Math.,
\textbf{478}
(1996),
127--176.


\bibitem[Gr] {Gr}M. Greenberg, {\it Rational points in discrete valuation
rings}, Publ. Math. I.H.E.S., {\bf 31} (1996), 59--64.
\bibitem[Hi] {Hi}M. Hickel, {\it Fonction de Artin et germes de courbes
trac{\'e}es sur un germe d'espace analytique}, Amer. J. Math., {\bf 115} (1993),
1299-1334.
\bibitem[Ig] {Ig}J.-I. Igusa, {\it Lectures on forms of higher degree}, Tata
Institute of Fundamental Research Lectures on Mathematics and Physics, {\bf 59}
(1978).
\bibitem[Ko] {Ko}M. Kontsevich, Lecture at Orsay (December 7, 1995).
\bibitem[Le] {Le}M. Lejeune-Jalabert, {\it Courbes trac{\'e}es sur un germe
d'hypersurface}, Amer. J. Math., {\bf 112} (1990), 525--568.
\bibitem[LeRe] {LeRe}M. Lejeune-Jalabert, A. Reguera-L{\'o}pez,
\textit{Arcs and wedges on sandwiched surface singularities},
Amer. J. Math., {\bf 121} (1999), 1191--1213.
\bibitem[Loe] {Loe}F. Loeser, {\it Fonctions d'Igusa $p$-adiques et polyn{\^o}mes de
Bernstein}, Amer. J. Math., {\bf 110} (1988), 1--21.
\bibitem[Loo] {Loo}E. Looijenga, {\it Motivic Measures}, in S{\'e}minaire Bourbaki,
expos{\'e} 874, Mars 2000.
\bibitem[Na] {Na}J. Nash Jr., {\it Arc structure of singularities}, Duke Math.
J., {\bf 81} (1995), 31--38.
\bibitem[Pa] {Pa}J. Pas, {\it Uniform $p$-adic cell decomposition and local zeta
functions}, J. Reine Angew. Math., {\bf 399} (1989), 137--172.
\bibitem[Re] {Re}M. Reid, {\it La correspondance de McKay}, in S{\'e}minaire
Bourbaki, expos{\'e} 867, Novembre 1999, 20 p., available at math. AG/9911165 (en
anglais).
\bibitem[Sa1] {Sa1}M. Saito, {\it Modules de Hodge polarisables}, Publ. Res.
Inst. Math. Sci., {\bf 24} (1988), 849--995.
\bibitem[Sa2] {Sa2}M. Saito, {\it Mixed Hodge modules}, Publ.
Res. Inst. Math. Sci.,
{\bf
26} (1990), 221--333.
\bibitem[Sa3] {Sa3}M. Saito, {\it Mixed Hodge modules and applications} in
Proceedings of the International Congress of Mathematicians, Vol. I, II (Kyoto,
1990), Math. Soc. Japan, Tokyo (1991), 725--734.
\bibitem[Sa4] {Sa4}M. Saito, {\it Hodge filtration on vanishing cycles},
preprint,
May 1998.



\bibitem[SGA 7] {SGA 7}
Groupes de monodromie en g{\'e}o\-m{\'e}\-trie al\-g{\'e}\-bri\-que,
di\-ri\-g{\'e} par {\sc A. Gro\-then\-dieck,} avec la collaboration de
{\sc M. Raynaud} et {\sc D. S. Rim}
pour la partie I et par
{\sc P. Deligne}
et {\sc N. Katz} pour la partie II,
{\it Lecture Notes in Math.,}
vol. 288 and 340,
Springer-Verlag 1972 and 1973.



\bibitem[Se] {Se}J.-P. Serre, {\it Facteurs locaux des fonctions z{\^e}ta des
vari{\'e}tes alg{\'e}\-bri\-ques (d{\'e}\-fi\-nitions et conjectures)}, S{\'e}minaire
Delange-Pisot-Poitou,
1969/70, n$^{\circ}$ {\bf 19}.
\bibitem[SeTh]  {SeTh}M. Sebastiani and R. Thom, {\it Un r{\'e}sultat sur la
monodromie}, Invent. Math., {\bf 13} (1971), 90--96.
\bibitem[Sc] {Sc}A. Scholl, {\it Classical motives} in Motives (Seattle, Wash.,
1991), Proc. Sympos. Pure Math., {\bf 55}, Part 1, Amer. Math. Soc.,
Providence, 1994, 163--187.
\bibitem[St1] {St1}J. Steenbrink,
\textit{Mixed Hodge structures on the vanishing cohomology},
in Real and Complex Singularities, Sijthoff and Noordhoff, Alphen aan 
den Rijn,
1977, 525--563.
\bibitem[St2]  {St2}J. Steenbrink, {\it The spectrum of hypersurface
singularities}, in Th{\'e}orie de Hodge, Luminy 1987, Ast{\'e}risque, {\bf 179-180}
(1989), 163--184.
\bibitem[St3] {St3}J. Steenbrink, {\it A summary of mixed Hodge theory}, American
Mathematical Society, Proc. Symp. Pure Math. {\bf 55}, Pt. 1, (1994), 31--41.
\bibitem[Va] {Va}A. Varchenko, {\it Asymptotic Hodge structure in the vanishing
cohomology}, Math. USSR Izvestija, {\bf 18} (1982), 469--512.   
\bibitem[Ve1] {Ve1}W. Veys, {\it Poles of Igusa's local zeta function and
monodromy}, Bull. Soc. Math. France, {\bf 121} (1993), 545--598.
\bibitem[Ve2] {Ve2}W. Veys, {\it Holomorphy of local zeta functions for curves},
Math. Ann. {\bf 295} (1993), 631--641.
\bibitem[Ve3] {Ve3}W. Veys, {\it Embedded resolution of singularities and Igusa's
local zeta function}, Academiae Analecta, to appear (survey paper, 46p).
\bibitem[Ve4] {Ve4}W. Veys, \textit{The topological zeta function associated to a function on a normal surface
germ}, Topology \textbf{38} (1999), 439--456. 
\bibitem[Ve5] {Ve5}W. Veys, {\it Zeta functions and `Kontsevich invariants' on
singular varieties}, Preprint Univ. Leuven 24 p. (1999).

\end{thebibliography}
\end{document}